# Bivariate Exponentaited Generalized Weibull-Gompertz Distribution

**M. A. EL-Damcese[1], Abdelfattah Mustafa[2], and M. S. Eliwa[2]**

[1] Mathematics Department, Faculty of Science, Tanta University, Egypt.

[2] Mathematics Department, Faculty of Science, Mansoura University, Mansoura 35516, Egypt.

**Abstract**

In this paper, we introduce a bivariate exponentaited generalized Weibull-Gompertz distribution. The model introduced here is of Marshall-Olkin type. Several properties are studied such as bivariate probability density function and it is marginal, moments, maximum likelihood estimation, joint reversed (hazard) function and joint mean waiting time and it is marginal. A real data set is analyzed and it is observed that the present distribution can provide a better fit than some other very well-known distributions.

**Keywords**: *Exponentaited generalized Weibull-Gompertz distribution, Marshall-Olkin, joint reversed (hazard) function.*

## 1. Introduction

Recently EL-Damcese1et al. (2014) introduced a five-parameters to generalize the generalized Weibull-Gompertz distribution by exponentiating generalized Weibull-Gompertz distribution which generalizes a lot of distributions such that generalized Gompertz distribution, generalized Weibull-Gompertz distribution, exponential power distribution, generalized exponential distribution, Weibull extension model of Chen (2000), Weibull extension model of Xie (2002) and etc. The exponentiation introduced an extra shape parameter in the model, which may yield more flexibility in the shape of the probability density function and hazard function.

This paper introduces a bivariate exponentaited generalized Weibull-Gompertz distribution (BEGWGD) by using the method of Marshall and Olkin (1967). In fact, shock models are used in reliability to describe different applications. Shocks can refer for example to damage caused to biological organs by illness or environmental causes of damage acting on a technical system. See for example, Sarhan and Balakrishnan (2007) studied Marshall and Olkin bivariate exponential distribution, Al- Khedhairi and El-Gohary (2008) presented a new class of bivariate Gompertz distributions, Kundu and Gupta (2009) proposed the bivariate generalized exponential distribution, El-Sherpieny et al. (2013) introduced a new bivariate generalized Gompertz distribution and Kundu and Gupta (2013) introduced Marshall–Olkin bivariate Weibull distribution.

[*]Corresponding authors: M. S. Eliwa[2]
E-mail: mseliwa@yahoo.com

The random variable $X$ is said to be has Exponentiated generalized Weibull-Gompertz distribution (EGWGD) with parameters $a, b, c, d$ and $\alpha$ if it has the following cumulative distribution function (CDF), as follows:

$$F_X(x;a,b,c,d,\alpha) = \left[1 - e^{-ax^b(e^{cx^d}-1)}\right]^{\alpha}, \quad a,b,c,d,\alpha > 0, \tag{1}$$

where $b, \alpha$ and $d$ are shape parameters, $a$ is scale parameter and $c$ is an acceleration parameter. The probability density function (PDF) of EGWGD$(a, b, c, d, \alpha)$ is

$$f_X(x;a,b,c,d,\alpha) = ab\alpha x^{b-1} e^{-ax^b(e^{cx^d}-1)+cx^d} \left(1 + \frac{cd}{b}x^d - e^{-cx^d}\right)\left[1 - e^{-ax^b(e^{cx^d}-1)}\right]^{\alpha-1}. \tag{2}$$

The survival function for the EGWGD$(a, b, c, d, \alpha)$ can be obtained as follows.

$$R(x;a,b,c,d,\alpha) = 1 - \left[1 - e^{-ax^b(e^{cx^d}-1)}\right]^{\alpha}, \quad x > 0 \tag{3}$$

## 2. Bivariate Exponentiated Generalized Weibull-Gompertz Distribution

In this section, we define a bivariate exponentiated generalized Weibull-Gompertz distribution (BEGWGD). We start with the joint cumulative distribution function and then derive the corresponding joint probability density function.

### 2.1 The joint cumulative distribution function

Let $U_1, U_2$ and $U_3$ be mutually independent random variables with the following distribution $U_i \sim \text{EGWGD}(a, b, c, d, \alpha_i), i = 1,2,3$.

Define the random variables $X_1$ and $X_2$ as

$$X_i = \max(U_i, U_3), i = 1, 2. \tag{4}$$

It is evident that the random variables $X_1$ and $X_2$ in (4) are dependent because the common random variable $U_3$.

We now study the joint distribution of the random variables $X_1$ and $X_2$. The following lemma gives the joint cumulative function of $X_1$ and $X_2$.

**Lemma 1.** The joint cumulative distribution function of the random variables $X_1$ and $X_2$ is

$$F_{X_1,X_2}(x_1,x_2) = \left[1 - e^{-ax_1^b(e^{cx_1^d}-1)}\right]^{\alpha_1}\left[1 - e^{-ax_2^b(e^{cx_2^d}-1)}\right]^{\alpha_2}\left[1 - e^{-az^b(e^{cz^d}-1)}\right]^{\alpha_3}, \tag{5}$$

where $z = \min(x_1, x_2), \alpha_i > 0, i = 1, 2, 3$.

**Proof.** Since

$$F_{X_1,X_2}(x_1,x_2) = P[X_1 \leq x_1, X_2 \leq x_2] = P[\max(U_1, U_3) \leq x_1, \max(U_2, U_3) \leq x_2]$$

$$= P[U_1 \leq x_1]P[U_2 \leq x_2]P[U_3 \leq \min(x_1, x_2)]$$

$$= F_{U_1}(x_1;a,b,c,d,\alpha_1)F_{U_2}(x_2;a,b,c,d,\alpha_2)F_{U_3}(z;a,b,c,d,\alpha_3)$$

Substituting from Equation (1) into the above equation, we get Equation (5). Which completes the proof of the lemma.

### 2.2. The joint probability density function

The following theorem gives the joint probability density function of the BEGWGD.

**Theorem 1.** If the joint cumulative distribution function of $(X_1, X_2)$ is as in Equation (5), the joint probability density function of the random variables $X_1$ and $X_2$ is given by

$$f_{X_1,X_2}(x_1,x_2) = \begin{cases} f_1(x_1,x_2) & \text{if } x_1 < x_2, \\ f_2(x_1,x_2) & \text{if } x_1 > x_2, \\ f_3(x,x) & \text{if } x_1 = x_2 = x, \end{cases} \tag{6}$$

where

$$f_1(x_1,x_2) = ab(\alpha_1+\alpha_3)x_1^{b-1}e^{-ax_1^b(e^{cx_1^d}-1)+cx_1^d}\left(1+\frac{cd}{b}x_1^d-e^{-cx_1^d}\right)\left[1-e^{-ax_1^b(e^{cx_1^d}-1)}\right]^{\alpha_1+\alpha_3-1}$$

$$\times ab\alpha_2 x_2^{b-1}e^{-ax_2^b(e^{cx_2^d}-1)+cx_2^d}\left(1+\frac{cd}{b}x_2^d-e^{-cx_2^d}\right)\left[1-e^{-ax_2^b(e^{cx_2^d}-1)}\right]^{\alpha_2-1},$$

$$f_2(x_1,x_2) = ab\alpha_1 x_1^{b-1}e^{-ax_1^b(e^{cx_1^d}-1)+cx_1^d}\left(1+\frac{cd}{b}x_1^d-e^{-cx_1^d}\right)\left[1-e^{-ax_1^b(e^{cx_1^d}-1)}\right]^{\alpha_1-1}\times$$

$$ab(\alpha_2+\alpha_3)x_2^{b-1}e^{-ax_2^b(e^{cx_2^d}-1)+cx_2^d}\left(1+\frac{cd}{b}x_2^d-e^{-cx_2^d}\right)\left[1-e^{-ax_2^b(e^{cx_2^d}-1)}\right]^{\alpha_2+\alpha_3-1},$$

and

$$f_3(x,x) = ab\alpha_3 x^{b-1}e^{-ax^b(e^{cx^d}-1)+cx^d}\left(1+\frac{cd}{b}x^d-e^{-cx^d}\right)\left[1-e^{-ax^b(e^{cx^d}-1)}\right]^{\alpha_1+\alpha_2+\alpha_3-1}.$$

**Proof.** Let us first assume that $x_1 < x_2$. In this case $F_{X_1,X_2}(x_1,x_2)$ in Equation (5) becomes

$$F_{X_1,X_2}(x_1,x_2) = \left[1-e^{-ax_1^b(e^{cx_1^d}-1)}\right]^{\alpha_1+\alpha_3}\left[1-e^{-ax_2^b(e^{cx_2^d}-1)}\right]^{\alpha_2}.$$

Then upon differentiation, we obtain the expression of $f_{X_1,X_2}(x_1,x_2) = \frac{\partial^2 F_{X_1,X_2}(x_1,x_2)}{\partial x_1 \partial x_2}$ to be $f_1(x_1,x_2)$ given above. Similarly, we find the expression of $f_{X_1,X_2}(x_1,x_2)$ to be $f_2(x_1,x_2)$ when $x_1 > x_2$. But to get $f_3(x,x)$ we use the following formula

$$\int_0^\infty \int_0^{x_2} f_1(x_1,x_2)dx_1 dx_2 + \int_0^\infty \int_0^{x_1} f_2(x_1,x_2)dx_2 dx_1 + \int_0^\infty f_3(x,x)dx = 1. \tag{7}$$

one can verify that

$$I_1 = \int_0^\infty \int_0^{x_2} f_1(x_1,x_2)\,dx_1 dx_2$$

$$I_1 = \int_0^\infty \int_0^{x_2} \left(ab(\alpha_1+\alpha_3)x_1^{b-1}e^{-ax_1^b(e^{cx_1^d}-1)+cx_1^d}\left(1+\frac{cd}{b}x_1^d-e^{-cx_1^d}\right)\left[1-e^{-ax_1^b(e^{cx_1^d}-1)}\right]^{\alpha_1+\alpha_3-1}\right.$$

$$ab\alpha_2 x_2^{b-1} e^{-ax_2^b\left(e^{cx_2^d}-1\right)+cx_2^d}\left(1+\frac{cd}{b}x_2^d-e^{-cx_2^d}\right)\left[1-e^{-ax_2^b\left(e^{cx_2^d}-1\right)}\right]^{\alpha_2-1}\Biggr)dx_1 dx_2$$

$$=\int_0^\infty ab\alpha_2 x_2^{b-1} e^{-ax_2^b\left(e^{cx_2^d}-1\right)+cx_2^d}\left(1+\frac{cd}{b}x_2^d-e^{-cx_2^d}\right)\left[1-e^{-ax_2^b\left(e^{cx_2^d}-1\right)}\right]^{\alpha_1+\alpha_2+\alpha_3-1} dx_2. \quad (8)$$

and,

$$I_2=\int_0^\infty ab\alpha_1 x_1^{b-1} e^{-ax_1^b\left(e^{cx_1^d}-1\right)+cx_1^d}\left(1+\frac{cd}{b}x_1^d-e^{-cx_1^d}\right)\left[-e^{-ax_1^b\left(e^{cx_1^d}-1\right)}\right]^{\alpha_1+\alpha_2+\alpha_3-1} dx_1. \quad (9)$$

Substituting from Equation (8) and Equation (9) into Equation (7), we have

$$\int_0^\infty f_3(x,x)dx=$$

$$=\int_0^\infty ab(\alpha_1+\alpha_2+\alpha_3)x^{b-1} e^{-ax^b\left(e^{cx^d}-1\right)+cx^d}\left(1+\frac{cd}{b}x^d-e^{-cx^d}\right)\left[1-e^{-ax^b\left(e^{cx^d}-1\right)}\right]^{\alpha_1+\alpha_2+\alpha_3-1} dx$$

$$-\int_0^\infty ab\alpha_2 x^{b-1} e^{-ax^b\left(e^{cx^d}-1\right)+cx^d}\left(1+\frac{cd}{b}x^d-e^{-cx^d}\right)\left[1-e^{-ax^b\left(e^{cx^d}-1\right)}\right]^{\alpha_1+\alpha_2+\alpha_3-1} dx$$

$$-\int_0^\infty ab\alpha_1 x^{b-1} e^{-ax^b\left(e^{cx^d}-1\right)+cx^d}\left(1+\frac{cd}{b}x^d-e^{-cx^d}\right)\left[1-e^{-ax^b\left(e^{cx^d}-1\right)}\right]^{\alpha_1+\alpha_2+\alpha_3-1} dx,$$

which readily yields

$$f_3(x,x)=ab\alpha_3 x^{b-1} e^{-ax^b\left(e^{cx^d}-1\right)+cx^d}\left(1+\frac{cd}{b}x^d-e^{-cx^d}\right)\left[1-e^{-ax^b\left(e^{cx^d}-1\right)}\right]^{\alpha_1+\alpha_2+\alpha_3-1}.$$

This completes the proof of the theorem.

### 3. The Marginal and Conditional Probability Density Function

In this section, we derive the marginal density function of $X_i$ and the conditional density function of $X_i|X_j, i\neq j=1,2$.

**Theorem 2.** The marginal probability density function of the random variables $X_i, i=1,2$ is given by

$$f_{X_i}(x_i)=ab(\alpha_i+\alpha_3)x_i^{b-1} e^{-ax_i^b\left(e^{cx_i^d}-1\right)+cx_i^d}\left(1+\frac{cd}{b}x_i^d-e^{-cx_i^d}\right)\left[1-e^{-ax_i^b\left(e^{cx_i^d}-1\right)}\right]^{\alpha_i+\alpha_3-1},$$

$$x_i>0, i=1,2. \quad (10)$$

**Proof.** Since

$$F_{X_i}(x_i)=P[X_i\leq x_i]=P[\max(U_i,U_3)\leq x_i]=\left[1-e^{-ax_i^b\left(e^{cx_i^d}-1\right)}\right]^{\alpha_i+\alpha_3} \quad (11)$$

Since $f_{X_i}(x_i)=\frac{d}{dx_i}F_{X_i}(x_i)$, then we get (10). Which completes the proof of the theorem.

**Remark 1.** From Equation (11), Theorem 2 and Equation (1), we note that $X_i$ has exponentaited generalized Weibull-Gompertz distribution with parameters $a,b,c,d,\alpha_i+\alpha_3$, $i=1,2$.

Having obtained the marginal probability density function of $X_1$ and $X_2$, we can now derive the conditional probability density functions as presented in the following theorem.

**Theorem 3.** The conditional probability density function of $X_i$ given that $X_j = x_j$, denoted by $f_{X_i|X_j}(x_i|x_j)$, $i \neq j = 1,2$, is given by

$$f_{X_i|X_j}(x_i|x_j) = \begin{cases} f^{(1)}_{X_i|X_j}(x_i|x_j), & \text{if } x_i < x_j, \\ f^{(2)}_{X_i|X_j}(x_i|x_j), & \text{if } x_i > x_j, \\ f^{(3)}_{X_i|X_j}(x_i|x_j), & \text{if } x_i = x_j, \end{cases} \quad (12)$$

where

$$f^{(1)}_{X_i|X_j}(x_i|x_j) = \frac{ab\alpha_j(\alpha_i + \alpha_3)x_i^{b-1}e^{-ax_i^b\left(e^{cx_i^d}-1\right)+cx_i^d}\left(1 + \frac{cd}{b}x_i^d - e^{-cx_i^d}\right)\left[1 - e^{-ax_i^b\left(e^{cx_i^d}-1\right)}\right]^{\alpha_i+\alpha_3-1}}{(\alpha_j + \alpha_3)\left[1 - e^{-ax_j^b\left(e^{cx_j^d}-1\right)}\right]^{\alpha_3}},$$

$$f^{(2)}_{X_i|X_j}(x_i|x_j) = ab\alpha_i x_i^{b-1}e^{-ax_i^b\left(e^{cx_i^d}-1\right)+cx_i^d}\left(1 + \frac{cd}{b}x_i^d - e^{-cx_i^d}\right)\left[1 - e^{-ax_i^b\left(e^{cx_i^d}-1\right)}\right]^{\alpha_i-1},$$

and

$$f^{(3)}_{X_i|X_j}(x_i|x_j) = \frac{\alpha_3\left[1 - e^{-ax^b\left(e^{cx^d}-1\right)}\right]^{\alpha_1-1}}{(\alpha_2 + \alpha_3)}.$$

**Proof.** The theorem follows readily upon substituting for the joint PDF of $(X_1, X_2)$ in Equation (6) and the marginal probability density function of $X_i, i = 1,2$, in Equation (10) in the following relation

$$f_{X_i|X_j}(x_i|x_j) = \frac{f_{X_i,X_j}(x_i,x_j)}{f_{X_j}(x_j)}, \quad i \neq j = 1,2.$$

## 4. The Mathematical Expectations

Based on the results presented in the last two sections, we can derive the mathematical expectations of $X_i, i = 1,2$.

**Theorem 4.** The *r-th* moments of $X_i$ is given by

$$E[X_i^r] = \frac{ab(\alpha_i + \alpha_3)}{d} \sum_{m=0}^{\infty}\sum_{j=0}^{\infty}\sum_{k=0}^{j}\sum_{l=0}^{\infty} \frac{(-1)^{m+j+k}c^l(a(1+m))^j}{j!\, l!\, d\, (ck)^{\frac{r+b(j+1)+ld}{d}}} \binom{j}{k}\binom{\alpha_i + \alpha_3 - 1}{m} \times$$

$$\left(((1+j)^l - j^l)\Gamma\left(\frac{r+b(j+1)+d(l-1)}{d} + 1\right) + \frac{d(1+j)^l}{kb}\Gamma\left(\frac{r+b(j+1)+ld}{d} + 1\right)\right). \quad (13)$$

**Proof:**

Substituting for $f_{X_i}(x_i)$ from (10), we get

$$E[X_i^r] = ab(\alpha_i + \alpha_3) \int_0^\infty x_i^{r+b-1} e^{-ax_i^b\left(e^{cx_i^d}-1\right)+cx_i^d} \left(1 + \frac{cd}{b} x_i^d - e^{-cx_i^d}\right) \times$$

$$\left[1 - e^{-ax_i^b(e^{cx_i^d}-1)}\right]^{\alpha_i+\alpha_3-1} dx_i$$

$$= ab(\alpha_i + \alpha_3)(I_1 - I_2 + \frac{cd}{b} I_3), \tag{14}$$

where

$$I_1 = \int_0^\infty x_i^{r+b-1} e^{-ax_i^b\left(e^{cx_i^d}-1\right)+cx_i^d} \left[1 - e^{-ax_i^b(e^{cx_i^d}-1)}\right]^{\alpha_i+\alpha_3-1} dx_i,$$

$$I_2 = \int_0^\infty x_i^{r+b-1} e^{-ax_i^b\left(e^{cx_i^d}-1\right)} \left[1 - e^{-ax_i^b(e^{cx_i^d}-1)}\right]^{\alpha_i+\alpha_3-1} dx_i,$$

$$I_3 = \int_0^\infty x_i^{r+b+d-1} e^{-ax_i^b\left(e^{cx_i^d}-1\right)+cx_i^d} \left[1 - e^{-ax_i^b(e^{cx_i^d}-1)}\right]^{\alpha_i+\alpha_3-1} dx_i.$$

Since $0 < \left[1 - e^{-ax_i^b(e^{cx_i^d}-1)}\right]^{\alpha_i+\alpha_3-1} < 1$ for $x_i > 0$, then by using the binomial series expansion we have

$$\left[1 - e^{-ax_i^b\left(e^{cx_i^d}-1\right)}\right]^{\alpha_i+\alpha_3-1} = \sum_{m=0}^\infty (-1)^m \binom{\alpha_i+\alpha_3-1}{m} e^{-max_i^b\left(e^{cx_i^d}-1\right)},$$

then

$$I_1 = \sum_{m=0}^\infty \sum_{j=0}^\infty \sum_{k=0}^j \sum_{l=0}^\infty \frac{(-1)^{m+j+k}(c(1+j))^l (a(1+m))^j}{j!\ l!} \binom{j}{k}\binom{\alpha_i+\alpha_3-1}{m} \times$$

$$\int_0^\infty x_i^{r+b+bj+ld-1} e^{-ckx_i^d} dx_i$$

$$= \sum_{m=0}^\infty \sum_{j=0}^\infty \sum_{k=0}^j \sum_{l=0}^\infty \frac{(-1)^{m+j+k}(c(1+j))^l (a(1+m))^j}{j!\ l!\ d\ (ck)^{\frac{r+b(j+1)+ld}{d}}} \binom{j}{k}\binom{\alpha_i+\alpha_3-1}{m} \times$$

$$\Gamma\left(\frac{r+b(j+1)+d(l-1)}{d}+1\right). \tag{15}$$

Similarly, we find that

$$I_2 = \sum_{m=0}^\infty \sum_{j=0}^\infty \sum_{k=0}^j \sum_{l=0}^\infty \frac{(-1)^{m+j+k}(cj)^l (a(1+m))^j}{j!\ l!\ d\ (ck)^{\frac{r+b(j+1)+ld}{d}}} \binom{j}{k}\binom{\alpha_i+\alpha_3-1}{m} \times$$

$$\Gamma\left(\frac{r+b(j+1)+d(l-1)}{d}+1\right). \tag{16}$$

$$I_3 = \sum_{m=0}^{\infty}\sum_{j=0}^{\infty}\sum_{k=0}^{j}\sum_{l=0}^{\infty} \frac{(-1)^{m+j+k}\bigl(c(1+j)\bigr)^l \bigl(a(1+m)\bigr)^j}{j!\,l!\,d\,(ck)^{\frac{r+b(j+1)+(l+1)d}{d}}} \binom{j}{k}\binom{\alpha_i+\alpha_3-1}{m} \times$$

$$\Gamma\left(\frac{r+b(j+1)+ld}{d}+1\right). \tag{17}$$

Substituting from Equation (15), Equation (16) and Equation (17) into Equation (14), we get Equation (13). This completes the proof.

## 5. Reliability Analysis

In this section we introduced some reliability measures, the joint reliability function, joint hazard rate, joint mean waiting time, joint reversed (hazard) function and it's marginal function.

### 5.1 Joint reliability function

In this subsection, we present the joint reliability function of $X_1$ and $X_2$, the reliability functions of $\min(X_1, X_2)$ and $\max(X_1, X_2)$. The following theorem gives the joint reliability function of the random variable $X_1$ and $X_2$.

**Theorem 5.** The joint reliability function of the random variables $X_1$ and $X_2$ is given by

$$R_{X_1,X_2}(x_1,x_2) = \begin{cases} R_1(x_1,x_2), & \text{if } x_1 < x_2 \\ R_2(x_1,x_2), & \text{if } x_2 < x_1 \\ R_3(x,x), & \text{if } x_1 = x_2 = x \end{cases} \tag{18}$$

where

$$R_1(x_1,x_2) = \left[1-e^{-ax_1^b\left(e^{cx_1^d}-1\right)}\right]^{\alpha_1+\alpha_3}\left[1-e^{-ax_2^b\left(e^{cx_2^d}-1\right)}\right]^{\alpha_2}\left(\left[1-e^{-ax_1^b\left(e^{cx_1^d}-1\right)}\right]^{-\alpha_1-\alpha_3}\times\right.$$

$$\left[1-e^{-ax_2^b\left(e^{cx_2^d}-1\right)}\right]^{-\alpha_2} - \left[1-e^{-ax_2^b\left(e^{cx_2^d}-1\right)}\right]^{-\alpha_2} - \left[1-e^{-ax_1^b\left(e^{cx_1^d}-1\right)}\right]^{-\alpha_1-\alpha_3}$$

$$\left.\left[1-e^{-ax_2^b\left(e^{cx_2^d}-1\right)}\right]^{\alpha_3}+1\right),$$

$$R_2(x_1,x_2) = \left[1-e^{-ax_1^b\left(e^{cx_1^d}-1\right)}\right]^{\alpha_1}\left[1-e^{-ax_2^b\left(e^{cx_2^d}-1\right)}\right]^{\alpha_2+\alpha_3}\left(\left[1-e^{-ax_1^b\left(e^{cx_1^d}-1\right)}\right]^{-\alpha_1}\times\right.$$

$$\left[1-e^{-ax_2^b\left(e^{cx_2^d}-1\right)}\right]^{-\alpha_2-\alpha_3} - \left[1-e^{-ax_1^b\left(e^{cx_1^d}-1\right)}\right]^{-\alpha_1} - \left[1-e^{-ax_1^b\left(e^{cx_1^d}-1\right)}\right]^{\alpha_3}\times$$

$$\left.\left[1-e^{-ax_2^b\left(e^{cx_2^d}-1\right)}\right]^{-\alpha_2-\alpha_3}+1\right),$$

and

$$R_3(x,x) = \left[1-e^{-ax^b\left(e^{cx^d}-1\right)}\right]^{\alpha_1+\alpha_2+\alpha_3}\times$$

$$\left(1+\left[1-e^{-ax^b\left(e^{cx^d}-1\right)}\right]^{-(\alpha_1+\alpha_2+\alpha_3)} - \left[1-e^{-ax^b\left(e^{cx^d}-1\right)}\right]^{-\alpha_2} - \left[1-e^{-ax^b\left(e^{cx^d}-1\right)}\right]^{-\alpha_1}\right).$$

**Proof**. The joint reliability function of $X_1$ and $X_2$ can be obtained by

$$R_{X_1,X_2}(x_1,x_2) = 1 - [F_{X_1}(x_1) + F_{X_2}(x_2) - F_{X_1,X_2}(x_1,x_2)] \quad (19)$$

substituting from Equation (10) in Equation (19), we get

$$R_{X_1,X_2}(x_1,x_2) = 1 - \left[1 - e^{-ax_1^b\left(e^{cx_1^d}-1\right)}\right]^{\alpha_1+\alpha_3} - \left[1 - e^{-ax_2^b\left(e^{cx_2^d}-1\right)}\right]^{\alpha_2+\alpha_3} +$$

$$\left[1 - e^{-ax_1^b\left(e^{cx_1^d}-1\right)}\right]^{\alpha_1} \left[1 - e^{-ax_2^b\left(e^{cx_2^d}-1\right)}\right]^{\alpha_2} \left[1 - e^{-az^b(e^{cz^d}-1)}\right]^{\alpha_3},$$

where $z = \min(x_1,x_2)$, $\alpha_i > 0$, $i = 1,2,3$, one can verify (18). This completes the proof.

In the mentioned applications $X_1$ and $X_2$ could be the lifetimes of two components, drought intensities for two regions, risks for two insurance events, exchange rates in two time periods and so on. So, it is important to know which of the two variables $X_1$ and $X_2$ is larger or smaller, so we want to obtain the distributions of $S = min(X_1, X_2)$ and $T = \max(X_1, X_2)$.

**Lemma 2**. The cumulative distribution functions of $S = \min(X_1, X_2)$ is given as

$$F_S(t) = \left[1 - e^{-at^b\left(e^{ct^d}-1\right)}\right]^{\alpha_1+\alpha_3} + \left[1 - e^{-at^b\left(e^{ct^d}-1\right)}\right]^{\alpha_2+\alpha_3} - \left[1 - e^{-at^b\left(e^{ct^d}-1\right)}\right]^{\alpha_1+\alpha_2+\alpha_3} \quad (20)$$

**Proof**. It is easy to prove that by using Equations (18) and (20).

**Lemma 3**. The cumulative distribution functions of $T = \max(X_1, X_2)$ is given as

$$F_T(t) = \left[1 - e^{-at^b\left(e^{ct^d}-1\right)}\right]^{\alpha_1+\alpha_2+\alpha_3} \quad (21)$$

**Proof**. Since

$$F_T(t) = P[T \leq t] = P[\max(X_1, X_2) \leq t] = P[X_1 \leq t, X_2 \leq t]$$
$$= p[U_1 \leq t]\, p[U_2 \leq t]\, p[U_3 \leq t] = \left[1 - e^{-at^b\left(e^{ct^d}-1\right)}\right]^{\alpha_1+\alpha_2+\alpha_3}.$$

### 5.2. Monotonicity of failure rate function

Let $(X_1, X_2)$ be two dimensional random variable with probability density function $f(x_1, x_2)$. Basu (1971) defined bivariate failure rate function (BVFR) $h(x_1, x_2)$ as

$$h(x_1, x_2) = \frac{f(x_1, x_2)}{R(x_1, x_2)}. \quad (22)$$

Then, by substituting from (6) and (18) in (22), we find that the bivariate hazard function when $x_1 < (>)x_2$ be

$$h(x_1, x_2) = \begin{cases} h_1(x_1, x_2), & \text{if } x_1 < x_2 \\ h_2(x_1, x_2), & \text{if } x_1 > x_2 \\ h_3(x_1, x_2), & \text{if } x_1 = x_2 = x \end{cases} \quad (23)$$

where

$$h_1(x_1, x_2) = \frac{f_1(x_1, x_2)}{R_1(x_1, x_2)}$$

$$= \left[(ab)^2 \alpha_2(\alpha_1 + \alpha_3)(x_1 x_2)^{b-1} e^{-ax_1^b\left(e^{cx_1^d}-1\right)+cx_1^d-ax_2^b\left(e^{cx_2^d}-1\right)+cx_2^d} \left(1 + \frac{cd}{b} x_1^d - e^{-cx_1^d}\right)\right.$$

$$\left. \times \left(1 + \frac{cd}{b} x_2^d - e^{-cx_2^d}\right)\right] \left[\rho_1 \left(1 - e^{-ax_1^b\left(e^{cx_1^d}-1\right)}\right)\left[\left(1 - e^{-ax_2^b\left(e^{cx_2^d}-1\right)}\right)\right]\right]^{-1},$$

$$h_2(x_1, x_2) = \frac{f_2(x_1, x_2)}{R_2(x_1, x_2)}$$

$$= \left[(ab)^2 \alpha_1(\alpha_2 + \alpha_3)(x_1 x_2)^{b-1} e^{-ax_1^b\left(e^{cx_1^d}-1\right)+cx_1^d-ax_2^b\left(e^{cx_2^d}-1\right)+cx_2^d} \left(1 + \frac{cd}{b} x_1^d - e^{-cx_1^d}\right)\right.$$

$$\left. \times \left(1 + \frac{cd}{b} x_2^d - e^{-cx_2^d}\right)\right] \left[\rho_2 \left(1 - e^{-ax_1^b\left(e^{cx_1^d}-1\right)}\right)\left(1 - e^{-ax_2^b\left(e^{cx_2^d}-1\right)}\right)\right]^{-1},$$

and

$$h_3(x,x) = \frac{f_3(x,x)}{R_3(x,x)}$$

$$= \frac{ab\alpha_3 x^{b-1} e^{-ax^b\left(e^{cx^d}-1\right)+cx^d}\left(1 + \frac{cd}{b} x^d - e^{-cx^d}\right)\left[1 - e^{-ax^b\left(e^{cx^d}-1\right)}\right]^{-1}}{\left(1 + \left[1 - e^{-ax^b\left(e^{cx^d}-1\right)}\right]^{-(\alpha_1+\alpha_2+\alpha_3)} - \left[1 - e^{-ax^b\left(e^{cx^d}-1\right)}\right]^{-\alpha_2} - \left[1 - e^{-ax^b\left(e^{cx^d}-1\right)}\right]^{-\alpha_1}\right)},$$

where

$$\rho_1 = \left\{\left[1 - e^{-ax_1^b\left(e^{cx_1^d}-1\right)}\right]^{-\alpha_1-\alpha_3}\left[1 - e^{-ax_2^b\left(e^{cx_2^d}-1\right)}\right]^{-\alpha_2} - \left[1 - e^{-ax_2^b\left(e^{cx_2^d}-1\right)}\right]^{-\alpha_2}\right.$$

$$\left. - \left[1 - e^{-ax_1^b\left(e^{cx_1^d}-1\right)}\right]^{-\alpha_1-\alpha_3}\left[1 - e^{-ax_2^b\left(e^{cx_2^d}-1\right)}\right]^{\alpha_3} + 1\right\},$$

$$\rho_2 = \left\{\left[1 - e^{-ax_1^b\left(e^{cx_1^d}-1\right)}\right]^{-\alpha_1}\left[1 - e^{-ax_2^b\left(e^{cx_2^d}-1\right)}\right]^{-\alpha_2-\alpha_3} - \left[1 - e^{-ax_1^b\left(e^{cx_1^d}-1\right)}\right]^{-\alpha_1}\right.$$

$$\left. - \left[1 - e^{-ax_1^b\left(e^{cx_1^d}-1\right)}\right]^{\alpha_3}\left[1 - e^{-ax_2^b\left(e^{cx_2^d}-1\right)}\right]^{-\alpha_2-\alpha_3} + 1\right\}.$$

Cox (1972) defined bivariate failure rate function as a vector which is useful to measure the total life span of a two component parallel system

$$h(\underline{x}) = (h(x), h_{12}(x_1|x_2), h_{21}(x_2|x_1)), \qquad (24)$$

where

$$h(x) = \frac{f_X(x)}{R(x)}\bigg|_{X=\min(X_1, X_2)},$$

$$h_{12}(x_1|x_2) = -\frac{\frac{\partial^2}{\partial x_1 \partial x_2} R(x_1, x_2)}{\frac{\partial}{\partial x_2} R(x_1, x_2)}, \qquad x_1 > x_2,$$

$$h_{21}(x_2|x_1) = -\frac{\frac{\partial^2}{\partial x_1 \partial x_2} R(x_1, x_2)}{\frac{\partial}{\partial x_1} R(x_1, x_2)}, \qquad x_1 < x_2.$$

The first element in the vector gives the failure function of the system using the information that both the component has survived beyond $x$. The second element gives the failure function span of the first component given that it has survived to an age $x_1$ and the other has failed at $x_2$. Similar argument holds for the third element.

**Lemma 4.** If $(X_1, X_2)$ is a BEGWG random vector, then the bivariate failure function as a vector in Equation (24) is $h(\underline{x})$ where

$$h(x) = \frac{1}{1 - F_S(x)} \frac{\partial}{\partial x} F_S(x),$$

$$h_{12}(x_1|x_2) = \frac{ab\alpha_1 x_1^{b-1} e^{-ax_1^b\left(e^{cx_1^d}-1\right)+cx_1^d}\left(1 + \frac{cd}{b}x_1^d - e^{-cx_1^d}\right)\left[1 - e^{-ax_1^b\left(e^{cx_1^d}-1\right)}\right]^{\alpha_1-1}}{1 - \left[1 - e^{-ax_1^b\left(e^{cx_1^d}-1\right)}\right]^{\alpha_1}},$$

$$h_{21}(x_2|x_1) = \frac{ab\alpha_2 x_2^{b-1} e^{-ax_2^b\left(e^{cx_2^d}-1\right)+cx_2^d}\left(1 + \frac{cd}{b}x_2^d - e^{-cx_2^d}\right)\left[1 - e^{-ax_2^b\left(e^{cx_2^d}-1\right)}\right]^{\alpha_2-1}}{1 - \left[1 - e^{-ax_2^b\left(e^{cx_2^d}-1\right)}\right]^{\alpha_2}},$$

**Proof.** It is easy to prove this Lemma by substituting from Equations (6), (18) and (20) into Equation (24). This completes the proof of the lemma.

Johnson and Kotz (1975) defined the hazard gradient as a vector $(h_1(x_1, x_2), h_2(x_1, x_2))$ where $h_1(x_1, x_2)$ is the hazard function of the conditional distribution of $X_1$ given $X_2 > x_2$ and $h_2(x_1, x_2)$ is the hazard function of the conditional distribution of $X_2$ given $X_1 > x_1$ as

$$\left.\begin{array}{l} h_1(x_1, x_2) = h(X_1|X_2 > x_2) = -\dfrac{\partial}{\partial x_1} \ell n\, R(x_1, x_2) \\[2mm] h_2(x_1, x_2) = h(X_2|X_1 > x_1) = -\dfrac{\partial}{\partial x_2} \ell n\, R(x_1, x_2) \end{array}\right\}. \quad (25)$$

Then by using Equations (18) and (25), we find the hazard gradient vector of BEGWGD when $x_1 < x_2$ is to be

$h_1(x_1, x_2) =$

$$\frac{\gamma_1(\alpha_1 + \alpha_3)\left(1 - \left[1 - e^{-ax_2^b\left(e^{cx_2^d}-1\right)}\right]^{\alpha_2}\right)\left[1 - e^{-ax_1^b\left(e^{cx_1^d}-1\right)}\right]^{\alpha_1+\alpha_3-1}}{1 - \left[1 - e^{-ax_1^b\left(e^{cx_1^d}-1\right)}\right]^{\alpha_1+\alpha_3} - \left[1 - e^{-ax_2^b\left(e^{cx_2^d}-1\right)}\right]^{\alpha_2+\alpha_3} + \left[1 - e^{-ax_1^b\left(e^{cx_1^d}-1\right)}\right]^{\alpha_1+\alpha_3}\left[1 - e^{-ax_2^b\left(e^{cx_2^d}-1\right)}\right]^{\alpha_2}},$$

And

$h_2(x_1, x_2) =$

$$\frac{\gamma_2\left((\alpha_2+\alpha_3)\left[1-e^{-ax_2^b\left(e^{cx_2^d}-1\right)}\right]^{\alpha_3} - \alpha_2\left[1-e^{-ax_1^b\left(e^{cx_1^d}-1\right)}\right]^{\alpha_1+\alpha_3}\right)\left[1-e^{-ax_2^b\left(e^{cx_2^d}-1\right)}\right]^{\alpha_2-1}}{1-\left[1-e^{-ax_1^b\left(e^{cx_1^d}-1\right)}\right]^{\alpha_1+\alpha_3} - \left[1-e^{-ax_2^b\left(e^{cx_2^d}-1\right)}\right]^{\alpha_2+\alpha_3} + \left[1-e^{-ax_1^b\left(e^{cx_1^d}-1\right)}\right]^{\alpha_1+\alpha_3}\left[1-e^{-ax_2^b\left(e^{cx_2^d}-1\right)}\right]^{\alpha_2}},$$

Similarly, we can find the hazard gradient vector of EGWG distribution when $x_1 > x_2$, so we can write the hazard gradient vector in the form

$h_1(x_1, x_2) =$

$$\frac{\gamma_1\left((\alpha_1+\alpha_3)\left[1-e^{-ax_1^b\left(e^{cx_1^d}-1\right)}\right]^{\alpha_3} - \alpha_1\left[1-e^{-ax_2^b\left(e^{cx_2^d}-1\right)}\right]^{\alpha_2+\alpha_3}\right)\left[1-e^{-ax_1^b\left(e^{cx_1^d}-1\right)}\right]^{\alpha_1-1}}{1-\left[1-e^{-ax_1^b\left(e^{cx_1^d}-1\right)}\right]^{\alpha_1+\alpha_3} - \left[1-e^{-ax_2^b\left(e^{cx_2^d}-1\right)}\right]^{\alpha_2+\alpha_3} + \left[1-e^{-ax_1^b\left(e^{cx_1^d}-1\right)}\right]^{\alpha_1}\left[1-e^{-ax_2^b\left(e^{cx_2^d}-1\right)}\right]^{\alpha_2+\alpha_3}},$$

and

$h_2(x_1, x_2) =$

$$\frac{\gamma_2(\alpha_2+\alpha_3)\left(1-\left[1-e^{-ax_1^b\left(e^{cx_1^d}-1\right)}\right]^{\alpha_1}\right)\left[1-e^{-ax_2^b\left(e^{cx_2^d}-1\right)}\right]^{\alpha_2+\alpha_3-1}}{1-\left[1-e^{-ax_1^b\left(e^{cx_1^d}-1\right)}\right]^{\alpha_1+\alpha_3} - \left[1-e^{-ax_2^b\left(e^{cx_2^d}-1\right)}\right]^{\alpha_2+\alpha_3} + \left[1-e^{-ax_1^b\left(e^{cx_1^d}-1\right)}\right]^{\alpha_1}\left[1-e^{-ax_2^b\left(e^{cx_2^d}-1\right)}\right]^{\alpha_2+\alpha_3}},$$

where

$$\gamma_i = abx_i^{b-1}e^{-ax_i^b\left(e^{cx_i^d}-1\right)+cx_i^d}\left(1+\frac{cd}{b}x_i^d - e^{-cx_i^d}\right), \quad i = 1,2.$$

### 5.3. The Mean waiting time

The reversed hazard function is closely related to another important random variable the waiting time. Indeed, as a condition of a failure in $[0, t]$ is already imposed while defining the reversed hazard function, it is of interest in different applications (actuarial science, reliability analysis) to describe the time, which had elapsed since the failure. The observations of waiting times, for instance, can be used for prediction (estimation) of the governing distribution function $F(t)$ will be shown later that the mean of the waiting time as a function of $t$ under certain assumptions defines uniquely $F(t)$. The waiting time could be of interest while describing different maintenance strategies. The mean waiting time function $m_w(t)$ for the marginal distribution of $X_1$ and $X_2$ can be denoted by

$$m_{w_i}(t) = \frac{1}{F_{X_i}(t)}\int_0^t F_{X_i}(x_i)dx_i \quad ; \quad i = 1,2, \tag{26}$$

where $F_{X_i}(t); i = 1,2$ be the cumulative distribution function to random variables $X_1$ and $X_2$ respectively so, the mean waiting time to join $(X_1, X_2)$ is

$$m_w(t_1, t_2) = \frac{1}{F(t_1, t_2)} \int_0^{t_1} \int_0^{t_2} F(x_1, x_2) \, dx_2 \, dx_1. \qquad (27)$$

The following Lemma obtains the mean waiting time to join $(X_1, X_2)$.

**Lemma 5.** The joint mean waiting time $m_w(t_1, t_2)$ to the random variables $X_1$ and $X_2$ is

$$m_w(t_1, t_2) = \begin{cases} m_{w_1}(t_1, t_2) & \text{if } t_1 < t_2, \\ m_{w_2}(t_1, t_2) & \text{if } t_1 > t_2, \\ m_{w_3}(t, t) & \text{if } t_1 = t_2 = t, \end{cases} \qquad (28)$$

where

$$m_{w_i}(t_1, t_2) = \frac{1}{F(t_1, t_2)} \sum_{m=0}^{\infty} \sum_{j=0}^{\infty} \sum_{k=0}^{j} \sum_{l=0}^{\infty} \frac{(am)^{2j}(c(j-k))^{2l}}{(j!\,l!)^2 (bj+dl+1)^2} \binom{j}{k}^2 \binom{\alpha_{3-i}}{m} \binom{\alpha_i + \alpha_3}{m} (t_1 t_2)^{bj+dl+1},$$

$i = 1,2$. And

$$m_{w_3}(t, t) = \frac{1}{F(t, t)} \sum_{m=0}^{\infty} \sum_{j=0}^{\infty} \sum_{k=0}^{j} \sum_{l=0}^{\infty} \frac{(-1)^{m+k+j}(am)^j (c(j-k))^l}{j!\,l!\,(bj+dl+1)} \binom{j}{k} \binom{\alpha_1 + \alpha_2 + \alpha_3}{m} t^{bj+dl+1}.$$

**Proof.** By using Equations (5) and (28) when $x_1 < x_2$, we get

$$m_{w_1}(t_1, t_2) = \frac{1}{F(t_1, t_2)} \int_0^{t_1} \int_0^{t_2} \left[1 - e^{-ax_1^b \left(e^{cx_1^d} - 1\right)}\right]^{\alpha_1 + \alpha_3} \left[1 - e^{-ax_2^b \left(e^{cx_2^d} - 1\right)}\right]^{\alpha_2} dx_2 \, dx_1$$

$$= \frac{1}{F(t_1, t_2)} \left\{ \left( \int_0^{t_1} \left[1 - e^{-ax_1^b \left(e^{cx_1^d} - 1\right)}\right]^{\alpha_1 + \alpha_3} dx_1 \right) \left( \int_0^{t_2} \left[1 - e^{-ax_2^b \left(e^{cx_2^d} - 1\right)}\right]^{\alpha_2} dx_2 \right) \right\},$$

but

$$\int_0^{t_2} \left[1 - e^{-ax_2^b \left(e^{cx_2^d} - 1\right)}\right]^{\alpha_2} dx_2 = \sum_{m=0}^{\infty} \sum_{j=0}^{\infty} \sum_{k=0}^{j} \sum_{l=0}^{\infty} \frac{(-1)^{m+k+j}(am)^j (c(j-k))^l}{j!\,l!\,(bj+dl+1)} \binom{j}{k} \binom{\alpha_2}{m} t_2^{bj+dl+1},$$

$$\int_0^{t_1} \left[1 - e^{-ax_1^b \left(e^{cx_1^d} - 1\right)}\right]^{\alpha_1 + \alpha_3} dx_1 = \sum_{m=0}^{\infty} \sum_{j=0}^{\infty} \sum_{k=0}^{j} \sum_{l=0}^{\infty} \frac{(-1)^{m+k+j}(am)^j (c(j-k))^l}{j!\,l!\,(bj+dl+1)} \binom{j}{k} \binom{\alpha_1 + \alpha_3}{m} t_1^{bj+dl+1}.$$

Therefore

$$m_{w_1}(t_1, t_2) = \frac{1}{F(t_1, t_2)} \sum_{m=0}^{\infty} \sum_{j=0}^{\infty} \sum_{k=0}^{j} \sum_{l=0}^{\infty} \frac{(am)^{2j}(c(j-k))^{2l}}{(j!\,l!)^2 (bj+dl+1)^2} \binom{j}{k}^2 \binom{\alpha_2}{m} \binom{\alpha_1 + \alpha_3}{m} (t_1 t_2)^{bj+dl+1}$$

Similarly, when $x_1 < x_2$ then

$$m_{w_2}(t_1, t_2) = \frac{1}{F(t_1, t_2)} \sum_{m=0}^{\infty} \sum_{j=0}^{\infty} \sum_{k=0}^{j} \sum_{l=0}^{\infty} \frac{(am)^{2j}(c(j-k))^{2l}}{(j!\,l!)^2 (bj+dl+1)^2} \binom{j}{k}^2 \binom{\alpha_1}{m} \binom{\alpha_2 + \alpha_3}{m} (t_1 t_2)^{bj+dl+1}.$$

Also; when $x_1 = x_2 = x$, we get

$$m_{w_3}(t,t) = \frac{1}{F(t,t)} \sum_{m=0}^{\infty} \sum_{j=0}^{\infty} \sum_{k=0}^{j} \sum_{l=0}^{\infty} \frac{(-1)^{m+k+j}(am)^j(c(j-k))^l}{j!\, l!\, (bj+dl+1)} \binom{j}{k} \binom{\alpha_1+\alpha_2+\alpha_3}{m} t^{bj+dl+1}.$$

Which complete the proof of this Lemma.

The following lemma obtains that we can get the marginal distribution of $X_1$ and $X_2$ by the marginal mean waiting time to $X_1$ and $X_2$ respectively. Let $F_{X_1}(t)$ and $F_{X_2}(t)$ be the marginal distribution of $X_1$ and $X_2$ respectively, and let $m_{w_1}(t)$ and $m_{w_2}(t)$ be the mean waiting time to marginal $X_1$ and $X_2$ respectively.

**Lemma 6.** Let $X_1$ and $X_2$ be two random variables have EGWG distribution, then the marginal distribution to $X_1$ and $X_2$ can be written in a form

$$F_{X_i}(t) = \frac{1}{m_{w_i}(t)} \sum_{m=0}^{\infty} \sum_{j=0}^{\infty} \sum_{k=0}^{j} \sum_{l=0}^{\infty} \frac{(-1)^{m+k+j}(am)^j(c(j-k))^l}{j!\, l!\, (bj+dl+1)} \binom{j}{k} \binom{\alpha_i+\alpha_3}{m} t^{bj+dl+1}, i=1,2. \quad (29)$$

**Proof:** Since

$$m_{w_i}(t) = \frac{1}{F_{X_i}(t)} \int_0^t F_{X_i}(x_i) dx_i = \frac{1}{F_{X_i}(t)} \int_0^t \left[1 - e^{-ax_i^b\left(e^{cx_i^d}-1\right)}\right]^{\alpha_i+\alpha_3} dx_i, \quad i = 1,2. \quad (30)$$

But by using the binomial expansion

$$\left[1 - e^{-ax_i^b(e^{cx_i^d}-1)}\right]^{\alpha_i+\alpha_3} = \sum_{m=0}^{\infty} (-1)^m \binom{\alpha_i+\alpha_3}{m} e^{-max_i^b(e^{cx_i^d}-1)},$$

then

$$m_{w_i}(t) = \frac{1}{F_{X_i}(t)} \sum_{m=0}^{\infty} \sum_{j=0}^{\infty} \sum_{k=0}^{j} \sum_{l=0}^{\infty} \frac{(-1)^{m+k+j}(am)^j(c(j-k))^l}{j!\, l!} \binom{j}{k} \binom{\alpha_i+\alpha_3}{m} \int_0^t x_i^{bj+dl} dx_i$$

$$= \frac{1}{F_{X_i}(t)} \sum_{m=0}^{\infty} \sum_{j=0}^{\infty} \sum_{k=0}^{j} \sum_{l=0}^{\infty} \frac{(-1)^{m+k+j}(am)^j(c(j-k))^l}{j!\, l!\, (bj+dl+1)} \binom{j}{k} \binom{\alpha_i+\alpha_3}{m} t^{bj+dl+1},$$

so

$$F_{X_i}(t) = \frac{1}{m_{w_i}(t)} \sum_{m=0}^{\infty} \sum_{j=0}^{\infty} \sum_{k=0}^{j} \sum_{l=0}^{\infty} \frac{(-1)^{m+k+j}(am)^j(c(j-k))^l}{j!\, l!\, (bj+dl+1)} \binom{j}{k} \binom{\alpha_i+\alpha_3}{m} t^{bj+dl+1}.$$

This completes the proof.

### 5.4. Bivariate reversed hazard rate function

The bivariate reversed hazard function is defined as the ratio of the probability density function and the corresponding cumulative distribution function. There are a number of parallels between bivariate reversed hazard function and bivariate hazard function. Both functions are nonnegative values, and both functions can be applied to continuous random variables as well as to discrete probability distributions. Both functions are a general alternative framework for examining probability distributions For example, given either a bivariate hazard function or a bivariate reverse hazard function, there is a corresponding cumulative distribution and a density function. we can denote:

1. The bivariate reversed hazard function by

$$r(x_1, x_2) = \frac{f(x_1, x_2)}{F(x_1, x_2)}. \tag{31}$$

Then, by using Equations (5), (6) and (31), we find that the bivariate reversed hazard function to EGWGD $r(x_1, x_2)$ be

$$r(x_1, x_2) = \begin{cases} r_1(x_1, x_2) & \text{if } x_1 < x_2, \\ r_2(x_1, x_2) & \text{if } x_1 > x_2, \\ r_3(x, x) & \text{if } x_1 = x_2 = x. \end{cases} \tag{32}$$

where

$$r_i(x_1, x_2) = (ab)^2 \alpha_{3-i}(\alpha_i + \alpha_3)(x_1 x_2)^{b-1} e^{-ax_1^b\left(e^{cx_1^d}-1\right)+cx_1^d-ax_2^b\left(e^{cx_2^d}-1\right)+cx_2^d}$$
$$\left[1 - e^{-ax_1^b\left(e^{cx_1^d}-1\right)}\right]^{-1} \left(1 + \frac{cd}{b}x_1^d - e^{-cx_1^d}\right)\left(1 + \frac{cd}{b}x_2^d - e^{-cx_2^d}\right)\left[1 - e^{-ax_2^b\left(e^{cx_2^d}-1\right)}\right]^{-1}, \quad i = 1,2.$$

and

$$r_3(x, x) = a\alpha_3 b x^{b-1} e^{-ax^b\left(e^{cx^d}-1\right)+cx^d} \left(1 + \frac{cd}{b}x^d - e^{-cx^d}\right)\left[1 - e^{-ax^b\left(e^{cx^d}-1\right)}\right]^{-1}.$$

2. The gradient vector of bivariate reversed hazard function by $r(x_1, x_2) = \left(r_1(x_1), r_2(x_2)\right)$, where

$$r_i(x_i) = \frac{f_{X_i}(x_i)}{F_{X_i}(x_i)} = \frac{\partial}{\partial x_i} \log F_{X_i}(x_i) \, ; \, i = 1,2. \tag{33}$$

Let $X_1$ and $X_2$ are two random variables having EGWG distribution, then the gradient vector of bivariate reversed hazard function can be denoted by

$$r_i(x_i) = ab(\alpha_i + \alpha_3)x_i^{b-1}e^{-ax_i^b\left(e^{cx_i^d}-1\right)+cx_i^d}\left(1 + \frac{cd}{b}x_i^d - e^{-cx_i^d}\right)\left[1 - e^{-ax_i^b\left(e^{cx_i^d}-1\right)}\right]^{-1}, i = 1,2.$$

## 6. Parameters Estimation

To estimate the unknown parameters of the BEGWG distribution, we use the method of maximum likelihood. Consider constant values to the parameters $a, b, c, d$ so, we want to estimate the other parameters $\alpha_1, \alpha_2$ and $\alpha_3$. Suppose $\left((X_{11}, X_{21}), (X_{12}, X_{22}), \dots, (X_{1n}, X_{2n})\right)$ is a random sample from BEGWG distribution where

$$n_1 = (i, X_{1i} < X_{2i}), \quad n_2 = (i, X_{1i} > X_{2i}), \quad n_3 = (i, X_{1i} = X_{2i} = X_i), \quad n = \sum_{j=1}^{3} n_j. \tag{34}$$

By using the equations (6) and (34), we find that the likelihood of the sample is given by

Bivariate Exponentaited Generalized Weibull-Gompertz Distribution.

$$l(\alpha_1, \alpha_2, \alpha_3) = \prod_{i=1}^{n_1} f_1(x_{1i}, x_{2i}) \prod_{i=1}^{n_2} f_2(x_{1i}, x_{2i}) \prod_{i=1}^{n_3} f_3(x_i, x_i). \tag{35}$$

The log-likelihood function can be written as

$$L(\alpha_1, \alpha_2, \alpha_3) = \sum_{i=1}^{n_1} \ln f_1(x_{1i}, x_{2i}) + \sum_{i=1}^{n_2} \ln f_2(x_{1i}, x_{2i}) + \sum_{i=1}^{n_3} \ln f_3(x_i, x_i)$$

$$= n_1 \ln(\alpha_2) + n_1 \ln(\alpha_1 + \alpha_3) + n_2 \ln(\alpha_2 + \alpha_3) + n_2 \ln(\alpha_1) + n_3 \ln(\alpha_3)$$

$$+ (\alpha_1 + \alpha_3 - 1) \sum_{i=1}^{n_1} \ln\left[1 - e^{-ax_{1i}^b\left(e^{cx_{1i}^d}-1\right)}\right] + (\alpha_2 - 1) \sum_{i=1}^{n_1} \ln\left[1 - e^{-ax_{2i}^b\left(e^{cx_{2i}^d}-1\right)}\right]$$

$$+ (\alpha_1 - 1) \sum_{i=1}^{n_2} \ln\left[1 - e^{-ax_{1i}^b\left(e^{cx_{1i}^d}-1\right)}\right] + (\alpha_2 + \alpha_3 - 1) \sum_{i=1}^{n_2} \ln\left[1 - e^{-ax_{2i}^b\left(e^{cx_{2i}^d}-1\right)}\right]$$

$$+ (\alpha_1 + \alpha_2 + \alpha_3 - 1) \sum_{i=1}^{n_3} \ln\left[1 - e^{-ax_i^b\left(e^{cx_i^d}-1\right)}\right] + (2(n_1 + n_2) + n_3) \ln(ab)$$

$$+ \sum_{i=1}^{n_1} \left\{-ax_{1i}^b\left(e^{cx_{1i}^d}-1\right) + cx_{1i}^d - ax_{2i}^b\left(e^{cx_{2i}^d}-1\right) + cx_{2i}^d\right\}$$

$$+ \sum_{i=1}^{n_2} \left\{-ax_{1i}^b\left(e^{cx_{1i}^d}-1\right) + cx_{1i}^d - ax_{2i}^b\left(e^{cx_{2i}^d}-1\right) + cx_{2i}^d\right\}$$

$$+ \sum_{i=1}^{n_3} \left(-ax_i^b\left(e^{cx_i^d}-1\right) + cx_i^d\right) + (b-1)\left(\sum_{i=1}^{n_1} \ln(x_{1i}x_{2i}) + \sum_{i=1}^{n_2} \ln(x_{1i}x_{2i}) + \sum_{i=1}^{n_3} \ln(x_i)\right)$$

$$+ \sum_{i=1}^{n_1} \ln\left(1 + \frac{cd}{b}x_{1i}^d - e^{-cx_{1i}^d}\right) + \sum_{i=1}^{n_1} \ln\left(1 + \frac{cd}{b}x_{2i}^d - e^{-cx_{2i}^d}\right) + \sum_{i=1}^{n_2} \ln\left(1 + \frac{cd}{b}x_{1i}^d - e^{-cx_{1i}^d}\right)$$

$$+ \sum_{i=1}^{n_2} \ln\left(1 + \frac{cd}{b}x_{2i}^d - e^{-cx_{2i}^d}\right) + \sum_{i=1}^{n_3} \ln\left(1 + \frac{cd}{b}x_i^d - e^{-cx_i^d}\right). \tag{36}$$

The normal equations are

$$\frac{\partial L}{\partial \alpha_1} = \frac{n_1}{\hat{\alpha}_1 + \hat{\alpha}_3} + \frac{n_2}{\hat{\alpha}_1} + \sum_{i=1}^{n_1} \ln\left[1 - e^{-ax_{1i}^b\left(e^{cx_{1i}^d}-1\right)}\right]$$

$$+ \sum_{i=1}^{n_2} \ln\left[1 - e^{-ax_{1i}^b\left(e^{cx_{1i}^d}-1\right)}\right] + \sum_{i=1}^{n_3} \ln\left[1 - e^{-ax_i^b\left(e^{cx_i^d}-1\right)}\right] = 0, \tag{37}$$

$$\frac{\partial L}{\partial \alpha_2} = \frac{n_1}{\hat{\alpha}_2} + \frac{n_2}{\hat{\alpha}_2 + \hat{\alpha}_3} + \sum_{i=1}^{n_1} \ln\left[1 - e^{-ax_{2i}^b\left(e^{cx_{2i}^d}-1\right)}\right]$$

$$+ \sum_{i=1}^{n_2} \ln\left[1 - e^{-ax_{2i}^b\left(e^{cx_{2i}^d}-1\right)}\right] + \sum_{i=1}^{n_3} \ln\left[1 - e^{-ax_i^b\left(e^{cx_i^d}-1\right)}\right] = 0, \tag{38}$$

$$\frac{\partial L}{\partial \alpha_3} = \frac{n_1}{\hat{\alpha}_1 + \hat{\alpha}_3} + \frac{n_2}{\hat{\alpha}_2 + \hat{\alpha}_3} + \frac{n_3}{\hat{\alpha}_3} + \sum_{i=1}^{n_1} \ln\left[1 - e^{-ax_{1i}^b\left(e^{cx_{1i}^d}-1\right)}\right]$$

$$+ \sum_{i=1}^{n_2} \ln\left[1 - e^{-ax_{2i}^b\left(e^{cx_{2i}^d}-1\right)}\right] + \sum_{i=1}^{n_3} \ln\left[1 - e^{-ax_i^b\left(e^{cx_i^d}-1\right)}\right] = 0. \quad (39)$$

The normal equations do not have explicit solution and they have obtained it numerically. The MLE of $\theta$, say $\hat{\alpha}_3(\hat{\alpha}_1, \hat{\alpha}_2)$ can be obtained as

$$\hat{\alpha}_3(\hat{\alpha}_1, \hat{\alpha}_2) = \frac{n_3}{\beta + \gamma - \sum_{i=1}^{n_1} \ln\left[1 - e^{-ax_{1i}^b\left(e^{cx_{1i}^d}-1\right)}\right] - \sum_{i=1}^{n_2} \ln\left[1 - e^{-ax_{2i}^b\left(e^{cx_{2i}^d}-1\right)}\right] - \sum_{i=1}^{n_3} \ln\left[1 - e^{-ax_i^b\left(e^{cx_i^d}-1\right)}\right]}$$
(40)

where

$$\beta = \frac{n_2}{\hat{\alpha}_1} + \sum_{i=1}^{n_1} \ln\left[1 - e^{-ax_{1i}^b\left(e^{cx_{1i}^d}-1\right)}\right] + \sum_{i=1}^{n_2} \ln\left[1 - e^{-ax_{1i}^b\left(e^{cx_{1i}^d}-1\right)}\right] + \sum_{i=1}^{n_3} \ln\left[1 - e^{-ax_i^b\left(e^{cx_i^d}-1\right)}\right],$$

$$\gamma = \frac{n_1}{\hat{\alpha}_2} + \sum_{i=1}^{n_1} \ln\left[1 - e^{-ax_{2i}^b\left(e^{cx_{2i}^d}-1\right)}\right] + \sum_{i=1}^{n_2} \ln\left[1 - e^{-ax_{2i}^b\left(e^{cx_{2i}^d}-1\right)}\right] + \sum_{i=1}^{n_3} \ln\left[1 - e^{-ax_i^b\left(e^{cx_i^d}-1\right)}\right].$$

So, the MLEs of $\hat{\alpha}_1$ and $\hat{\alpha}_2$ can be obtained by solving two nonlinear Equations (37) and (38) by using Equation (40).

### 6.1. Asymptotic confidence bounds

In this section, we consider a constant value to the parameters $a, b, c$ and $d$ which take the values $0.1, 0.2, 0.2$ and $0.5$ respectively. Also, we derive the asymptotic confidence intervals of the parameters $\hat{\alpha}_i > 0, i = 1,2,3$ by using variance covariance matrix $I_0^{-1}$ see Lawless (2003), where $I_0^{-1}$ is the inverse of the observed information matrix

$$I_0^{-1} = -\begin{pmatrix} \frac{\partial^2 L}{\partial \alpha_1^2} & \frac{\partial^2 L}{\partial \alpha_1 \alpha_2} & \frac{\partial^2 L}{\partial \alpha_1 \alpha_3} \\ \frac{\partial^2 L}{\partial \alpha_2 \alpha_1} & \frac{\partial^2 L}{\partial \alpha_2^2} & \frac{\partial^2 L}{\partial \alpha_2 \alpha_3} \\ \frac{\partial^2 L}{\partial \alpha_3 \alpha_1} & \frac{\partial^2 L}{\partial \alpha_3 \alpha_2} & \frac{\partial^2 L}{\partial \alpha_3^2} \end{pmatrix}^{-1}, \quad (41)$$

thus

$$I_0^{-1} = \begin{bmatrix} var(\hat{\alpha}_1) & cov(\hat{\alpha}_1, \hat{\alpha}_2) & cov(\hat{\alpha}_1, \hat{\alpha}_3) \\ cov(\hat{\alpha}_2, \hat{\alpha}_1) & var(\hat{\alpha}_2) & cov(\hat{\alpha}_2, \hat{\alpha}_3) \\ cov(\hat{\alpha}_3, \hat{\alpha}_1) & cov(\hat{\alpha}_3, \hat{\alpha}_2) & var(\hat{\alpha}_3) \end{bmatrix}. \quad (42)$$

The derivatives in $I_0$ are given as follows:

$$\frac{\partial^2 L}{\partial \alpha_1^2} = -\left(\frac{n_1}{(\alpha_1 + \alpha_3)^2} + \frac{n_2}{\alpha_1^2}\right), \quad \frac{\partial^2 L}{\partial \alpha_1 \alpha_2} = 0, \quad L_{13} = \frac{\partial^2 L}{\partial \alpha_1 \alpha_3} = \frac{-n_1}{(\alpha_1 + \alpha_3)^2},$$

Bivariate Exponentaited Generalized Weibull-Gompertz Distribution.

$$\frac{\partial^2 L}{\partial \alpha_2^2} = -\left(\frac{n_2}{(\alpha_2+\alpha_3)^2}+\frac{n_1}{\alpha_2^2}\right), \quad \frac{\partial^2 L}{\partial \alpha_2 \alpha_3} = \frac{-n_2}{(\alpha_2+\alpha_3)^2},$$

$$\frac{\partial^2 L}{\partial \alpha_3^2} = -\left(\frac{n_1}{(\alpha_1+\alpha_3)^2}+\frac{n_2}{(\alpha_2+\alpha_3)^2}+\frac{n_3}{\alpha_3^2}\right).$$

We can derive the $(1-\delta)100\%$ confidence intervals of the parameters $\hat{\alpha}_i > 0, i = 1,2,3$ by using variance covariance matrix as in the following forms

$$\hat{\alpha}_i \pm Z_{\frac{\delta}{2}}\sqrt{var(\hat{\alpha}_i)}, \quad i = 1,2,3.$$

where $Z_{\delta/2}$ is the upper $(\delta/2)$th percentile of the standard normal distribution.

### 6.2 Data analysis

In this section, we present the analysis of a real data set using the BEGWG model and compare it with the other fitted model like bivariate exponentiated Gompertez distribution BEGD($\alpha_1, \alpha_2, \alpha_3, \lambda$). The following data represent the American Football (National Football League) League data and they are obtained from the matches played on three consecutive weekends in 1986. This data is available in Csorgo and Welsh (1989). It is consider a bivariate data set where $X_1$ represents the "game time" to the first points scored by kicking the ball between goal posts. $X_2$ represents the "game time" to the first points scored by moving the ball into the end zone. These times are of interest to a casual spectator who wants to know how long one has to wait to watch a touchdown or to a spectator who is interested only at the beginning stages of a game. The variables have the following structure:

1. $X_1 < X_2$ : Means that the first score is a field goal.
2. $X_1 > X_2$ : Means the first score is an unconverted safety or touchdown.
3. $X_1 = X_2$ : Means the first score is a converted touchdown.

The data (scoring times in minutes and seconds) are represented in the following Table 1.

Table 1. American Football (National Football League) League data.

| $X_1$ | $X_2$ | $X_1$ | $X_2$ | $X_1$ | $X_2$ |
|---|---|---|---|---|---|
| 2.05 | 3.98 | 5.78 | 25.98 | 10.40 | 10.25 |
| 9.05 | 9.05 | 13.80 | 49.75 | 2.98 | 2.98 |
| 0.85 | 0.85 | 7.25 | 7.25 | 3.88 | 6.43 |
| 3.43 | 3.43 | 4.25 | 4.25 | 0.75 | 0.75 |
| 7.78 | 7.78 | 1.65 | 1.65 | 11.63 | 17.37 |
| 10.57 | 14.28 | 6.42 | 15.08 | 1.38 | 1.38 |
| 7.05 | 7.05 | 4.22 | 9.48 | 10.53 | 10.53 |
| 2.58 | 2.58 | 15.53 | 15.53 | 12.13 | 12.13 |
| 7.23 | 9.68 | 2.90 | 2.90 | 14.58 | 14.58 |
| 6.85 | 34.58 | 7.02 | 7.02 | 11.82 | 11.82 |
| 32.45 | 42.35 | 6.42 | 6.42 | 5.52 | 11.27 |
| 8.53 | 14.57 | 8.98 | 8.98 | 19.65 | 10.70 |
| 31.13 | 49.88 | 10.15 | 10.15 | 17.83 | 17.83 |
| 14.58 | 20.57 | 8.87 | 8.87 | 10.85 | 38.07 |

The BEGWGD model is used to fit this data set. The MLE(s) of the unknown parameter(s), the value of log – likelihood (L), Akaike information criterion (AIC), correct Akaike information criterion (CAIC) and Bayesian information criterion (BIC) test statistic two different models are given in Table 1

**Table 1. The MLE(s) of the parameter(s), L, AIC, CAIC and BIC**

| The Model | MLE(s) | - L | AIC | CAIC | BIC |
|---|---|---|---|---|---|
| BEGD $(\hat{\alpha}_1, \hat{\alpha}_2, \hat{\alpha}_3, \hat{\lambda})$ | $\hat{\alpha}_1 = 0.043, \hat{\alpha}_1 = 0.528$ $\hat{\alpha}_3 = 1.037, \hat{\lambda} = 0.787$ | 370.41 | 748.82 | 749.90 | 377.88 |
| BEGWGD $(\hat{\alpha}_1, \hat{\alpha}_2, \hat{\alpha}_3)$ | $\hat{\alpha}_1 = 0.032\,3, \hat{\alpha}_1 = 0.186$ $\hat{\alpha}_3 = 0.406$ | 354.03 | 714.06 | 714.69 | 359.63 |

Table 1 show that BEGWGD model is the best distribution because it has the smallest value of AIC, CAIC and BIC test. By substituting the MLE of unknown parameters in Equation (37), we get estimation of the variance covariance matrix as

$$I_0^{-1} = \begin{bmatrix} 0.0005173 & 0.0000052 & -0.000426 \\ 0.0000052 & 0.0021423 & -0.000125 \\ -0.000426 & -0.000125 & 0.0102564 \end{bmatrix}.$$

The approximate 95% two sided confidence interval of $\hat{\alpha}_1, \hat{\alpha}_2$ and $\hat{\alpha}_3$ are [0,0.077], [0.0955,0.277] and [0.207, 0.605] respectively.

**References**

[1] Al-Khedhairi, A. and El-Gohary, A. (2008). "A new class of bivariate Gompertz distributions" *Internatinal Journal of Mathematics Analysis*, 2(5), 235 – 253.

[2] Basu, A.P.(1971)."Bivariate failure rate". *American Statistics Association*, 66,103-104.

[3] Chen, Z. (2000). "A new two-parameter lifetime distribution with bathtub shape or increasing failure rate function". *Statistics and Probability Letters*, 49, 155–161.

[4] Cox, D.R. (1972). "Regression models and life tables". *Royal Statistics Society*, 34, 187-220.

[5] Csorgo, S. and Welsh, A. H. (1989). "Testing for exponential and Marshall-Olkin distribution". *Journal of Statistical Planning and Inference*, 23, 287-300.

[6] EL-Damcese, M.A., Mustafa, A. and Eliwa, M. S. (2014). " Exponentaited generalized Weibull-Gompertz distribution ". *http://arxiv.org/abs/1412.0705*.

[7] El-Gohary, A., Alshamrani, A. and Al-Otaibi, A. N. (2013). "The Generalized Gompertz Distribution". *Journal of Applied Mathematical Modelling*, 37(1-2), 13-24.

[8] El-Sherpieny, E. A., Ibrahim, S. A., and Bedar, R. E. (2013). "A new bivariate generalized Gompertz distribution". *Asian Journal of Applied Sciences*, 1-4, 2321 – 0893.

[9] Johnson, N.L. and Kotz, S. (1975). "A vector valued multivariate hazard rate". *Journal of Multivariate Analysis*, 5, 53-66.

[10] Kundu, D. and Gupta, R. D. (2009). "Bivariate generalized exponential distribution". *Journal of Multivariate Analysis*, 100(4), 581-593.

[11] Kundu, D., Gupta, K. (2013)." Bayes estimation for the Marshall–Olkin bivariate Weibull distribution". *Journal of Computational Statistics and Data Analysis*, 57(1), 271–281.

[12] Lawless, J. F. (2003). "Statistical Models and Methods for Lifetime Data". *John Wiley and Sons*, New York, 20, 1108-1113.

[13] Marshall, A. W. and Olkin, I. A. (1986). "A multivariate exponential distribution". *Journal of the American Statistical Association*, 62, 30-44.

[14] Sarhan, A. and Balakrishnan, N. (2007). "A new class of bivariate distributions and its mixture". Journal *of the Multivariate Analysis*, 98, 1508-1527.

[15] Xie, M., Tang, Y., and Goh, T. N. (2002). "A modified Weibull extension with bathtub-shaped failure rate function". *Reliability Engineering and System Safety*, 76, 279–285.